\documentclass[a4paper,11pt,twoside,reqno]{amsart}
\usepackage{amsmath,amsthm,amsfonts,amssymb,amscd,enumerate,verbatim,calc}

\date{}
\RequirePackage[dvips]{graphicx} \textheight 20cm
\begin{document}
\centerline{}
\centerline{\large \bf CR-SUBMANIFOLDS OF SOME LORENTZIAN}
\centerline{\large \bf MANIFOLDS AND K-MANIFOLDS}
\centerline{}
\centerline{PAYEL KARMAKAR}
\centerline{Department of Mathematics,}
\centerline{Jadavpur University, Kolkata-700032, India.}
\centerline{E-mail: payelkarmakar632@gmail.com}
\centerline{}
\newtheorem{Theorem}{\quad Theorem}[section]
\newtheorem{Definition}[Theorem]{\quad Definition}
\newtheorem{Corollary}[Theorem]{\quad Corollary}
\newtheorem{Lemma}[Theorem]{\quad Lemma}
\newtheorem{Example}[Theorem]{\emph{Example}}
\newtheorem{Proposition}[Theorem]{Proposition}
\centerline{\textbf{Abstract}}
In this paper I have studied about CR(Cauchy-Riemann)-submanifolds of Lorentzian Concircular Structure manifold ($(LCS)_n$-manifold), Lorentzian Para-Sasakian(LP)-cosymplectic manifold, S-manifold and Generalized Kenmotsu (GKM) manifold. I have discussed some results regarding distribution, structure vector field, totally geodesic submanifold, leaf etc.. I have obtained results on totally umbilical contact CR-submanifold where the anti-invariant distribution has some properties. Next, I have studied some results about $D$-totally geodesic CR-submanifold ($D$ is the distribution), a contact CR-submanifold, $D^\perp$-totally geodesic CR-submanifold, $\xi$-horizontal CR-submanifold where the distribution is integrable (here $\xi$ is the structure vector field). Also I have proved some results on $D$-umbilic CR-submanifold, mixed totally geodesic CR-submanifold, foliate $\xi$-horizontal mixed totally geodesic CR-submanifold, leaf of the distribution, totally geodesic leaf, CR-product etc.. I have given an example of a CR-submanifold of an $(LCS)_n$-manifold and at last, I have given an example of a GKM manifold.
\let\thefootnote\relax\\\\
\let\thefootnote\relax
$~~~~~~~~~~~~~~~~~~~~~~~~~~~~~~~~~$\section{\large Introduction}
\
\par{In 1930, Schouten and Van Dantzing tried to transfer the results of Differential Geometry of spaces with Riemannian metric and affine connection to the case of spaces with complex structure. Using the complex structure and differential 1-form on a manifold, a great deal of work is carried out on these manifolds from 1960 onwards. These are known as contact manifolds and are odd dimensional. One can obtain different structures like Sasakian, Kenmotsu, Lorentzian etc., by providing additional conditions to the contact structure.}\\\
\par{A (2n+1)-dimensional smooth manifold $(M,g)$ is called almost contact if it admits an almost contact structure $(\phi,\xi,\eta)$, where $\phi$ is a tensor field of type (1,1), $\xi$ is a vector field and $\eta$ is a 1-form satisfying}
\begin{equation}
\phi^{2}X=-X+\eta(X)\xi,~\eta \circ \phi=0,~\phi\xi=0,~\eta(\xi)=1 \tag{1.1}
\end{equation}
for all vector fields $X,Y$ on $M$.\\\
\par{A differentiable $C^{\infty}$ (2n+1)-dimensional manifold $M^{(2n+1)}$ is called a contact manifold if it carries a global differential 1-form $\eta$ such that}
\begin{equation}
\eta \wedge(d\eta)^{n}\neq0 \tag{1.2}
\end{equation}
everywhere on $M^{2n+1}$ where the exponent denotes the nth exterior power.\\\
\\\
\par{An n-dimensional Lorentzian manifold $M$ is a smooth connected paracompact Hausdorff manifold with a Lorentzian metric $g$, i.e., $M$ admits a smooth symmetric tensor field $g$ of type (0,2) such that for each point the tensor $g_p:T_pM \times  T_pM\rightarrow \mathbb{R}$ is a non-degenerate inner-product of signature (-,+,...,+), $T_pM$ denotes the tangent vector space $M$ at $p$ and $\mathbb{R}$ is the real number space. In other words, a matrix representation of $g_p$ has one eigen value negative and all other eigen values positive. A non-zero vector $X_p\in T_pM$ is known to be spacelike, null or lightlike, non-spacelike or timelike according as $g_p(X_p,X_p)>0,=0$ or $<0$ respectively.}\\\
\par{If $M$ is a differentiable manifold of dimension n, and there exists a $(\phi,\xi,\eta)$ structure satisfying}
\begin{equation}
\phi^2=I-\eta\otimes\xi,~\eta(\xi)=1,~\phi(\xi)=0,~\eta \circ \phi=0, \tag{1.3}
\end{equation}
then $M$ is called an almost paracontact manifold.\\\
\par{In an almost paracontact structure $(\phi,\xi,\eta,g)$,}
\begin{equation}
g(X,\phi Y)=g(\phi X,Y), \tag{1.4}
\end{equation}
\begin{equation}
2\phi(X,Y)=(\bar{\nabla}_X \eta)Y+(\bar{\nabla}_Y \eta)X, \tag{1.5}
\end{equation}
where $\phi$ is a tensor of type (1,1),$\xi$ is a vector field, $\eta$ is a 1-form and $g$ is Lorentzian metric satisfying
\begin{equation}
\phi^2X=X+\eta(X)\xi,~\eta \circ \phi=0,~\phi\xi=0,~\eta(\xi)=-1, \tag{1.6}
\end{equation}
\begin{equation}
g(\phi X,\phi Y)=g(X,Y) + \eta(X)\eta(Y), g(X,\xi)=\eta(X) \tag{1.7}
\end{equation}
for all vector fields $X,Y$ on $M$.\\\
\\\
\par{In a Lorentzian manifold $(M,g)$, a vector field $P$ defined by $g(X,P)=A(X)$ for any
$X\in \Gamma(TM)$, is called con-circular if}
\begin{equation}
(\bar{\nabla}_X A)(Y)=\alpha \{g(X,Y)+\omega(X)A(Y)\}, \tag{1.8}
\end{equation}
where $\alpha$ is a non-zero scalar and $\omega$ is a closed 1-form and $\bar{\nabla}$ denotes the operator of covariant differentiation of $M$ with respect to the Lorentzian metric $g$.\\\
\par{Let $M$ admits a unit timelike concircular vector field $\xi$, called the structure vector field of the manifold, then $g(\xi,\xi)=-1$, since $\xi$ is a unit concircular vector field, it follows that $\exists$ a non-zero 1-form $\eta$ such that $g(X,\xi)=\eta(X)$. The following equations hold$-$}
\begin{equation}
(\bar{\nabla}_X \eta)Y=\alpha [g(X,Y) + \eta(X)\eta(Y)],~\alpha \neq 0, \tag{1.9}
\end{equation}
\begin{equation}
\bar{\nabla}_X \alpha= X_\alpha =d\alpha(X)=\rho\eta(X) \tag{1.10}
\end{equation}
for all vector fields $X,Y$ and $\alpha$ is a non-zero scalar function related to $\rho$ by $\rho=-(\xi\alpha)$.\\\ \par{Let $\phi X= \frac{1}{\alpha} \bar{\nabla}_X \xi$, from which it follows that $\phi$ is a symmetric (1,1) tensor and call it the structure tensor on the manifold. Thus the Lorentzian manifold $M$ together with unit timelike concircular vector field $\xi$, its associated 1-form $\eta$ and a (1,1) tensor field $\phi$ is called a Lorentzian Concircular Structure  manifold i.e, $(LCS)_n$-manifold. Specially if $\alpha=1$, then we obtain LP-Sasakian structure of Matsumoto[15]. In $(LCS)_n$-manifold $(n>2)$, the following relations hold$-$}
\begin{equation}
\phi^2 =I+\eta\otimes\xi,~\eta(\xi)=-1, \tag{1.11}
\end{equation}
where $I$ denotes the identity transformation of the tangent space $TM$. Also in an $(LCS)_n$-manifold, the following relations are satisfied $\forall$ $X,Y\in T(M)-$
\begin{equation}
\phi\xi=0,~\eta \circ \phi=0,~g(X,\phi Y)=g(\phi X,Y),~rank\phi=2n, \tag{1.12}
\end{equation}
\begin{equation}
g(\phi X,\phi Y)=g(X,Y)+\eta(X)\eta(Y), g(X,\xi)=\eta(X), \tag{1.13}
\end{equation}
\begin{equation}
\bar{R}(X,Y)\xi=(\alpha^2 -\rho)[\eta(Y)X-\eta(X)Y], \tag{1.14}
\end{equation}
\begin{equation}
\bar{R}(\xi,X)\xi=(\alpha^2 -\rho)[\eta(X)\xi+X]. \tag{1.15}
\end{equation}
\par{Also $(LCS)_n$-manifold satisfies$-$}
\begin{equation}
(\bar{\nabla}_X \phi)Y=\alpha [g(X,Y)\xi + 2\eta(X)\eta(Y)\xi + \eta(Y)X], \tag{1.16}
\end{equation}
\begin{equation}
\bar{R}(X,Y)Z= \phi R(X,Y)Z + (\alpha^2 -\rho)\{g(Y,Z)\eta(X)-g(X,Z)\eta(Y)\}, \tag{1.17}
\end{equation}
\begin{equation}
\bar{\nabla}_X \xi=\alpha \phi X. \tag{1.18}
\end{equation}
\par{A Lorentzian manifold $\bar{M}$ with the structure tensor $(\phi,\xi,\eta,g)$ satisfying (1.11), (1.12), (1.13) is called a Lorentzian almost paracontact manifold. In this case,}
\begin{equation}
(\bar{\nabla}_X\phi)Y=g(X,Y)\xi + 2\eta(X)\eta(Y)\xi, \tag{1.19}
\end{equation}
where $\bar{\nabla}$ is the covariant differentiation with respect to $g$.\\\
\par{A Lorentzian almost paracontact manifold is called an LP-cosymplectic manifold if}
\begin{equation}
(\bar{\nabla}_X \phi)Y=0, \bar{\nabla}_X \xi=0. \tag{1.20}
\end{equation}
\par{Tangent sphere bundle of any Riemannian manifold of constant sectional curvature with its standard metric structure is an LP-cosymplectic manifold.}\\\
\\\
\par{In 1963, Yano introduced notion of $\phi$-structure on a $C^{\infty}$ (2n+s) dimensional manifold $M$, as a non-vanishing tensor field $\phi$ of type (1,1) on $M$ which satisfies $\phi^3 +\phi=0$ and has constant rank $r=2n$. In 1970, Goldberg and Yano[5] defined globally framed $\phi$-structures, for which the sub-bundle $Ker \phi$ is parallelizable. Then $\exists$ a global frame ${\xi_1,...,\xi_s}$ for the sub-bundle $Ker\phi$, (the vector fields $\xi_1,...,\xi_s$ are called the structure vector fields) with dual 1-forms $\eta_1,...,\eta_s$ such that}
\begin{equation}
g(\phi X,\phi Y)=g(X,Y)-\sum_{\alpha=1}^{s} \eta_\alpha (X)\eta_\alpha (Y) \tag{1.21}
\end{equation}
for any vector fields $X,Y$ in $M$, then this structure is called a metric $\phi$-structure.\\\
\par{For S-manfolds, the following relations hold$-$}
\begin{equation}
\eta_\alpha (\xi_\beta)=\delta_{\alpha\beta},~\phi\xi_\alpha =0,~\eta_\alpha \circ \phi =0,~\tag{1.22}
\end{equation}
\begin{equation}
g(X,\phi Y)=-g(\phi X,Y), g(X,\xi_\alpha)=\eta_\alpha (X), \tag{1.23}
\end{equation}
\begin{equation}
(\bar{\nabla}_X \phi)Y= \sum_\alpha \{g(\phi X,\phi Y)\xi_\alpha +\eta_\alpha (Y)\phi^2 X\}, \tag{1.24}
\end{equation}
\begin{equation}
\bar{\nabla}_X \xi_\alpha =-\phi X \tag{1.25}
\end{equation}
$~~~~$ $\forall$ $X,Y\in TM,~\alpha,\beta=1,...,s$.\\\
\\\
\par{Let $M$ be (2n+s) dimensional K-manifold with structure tensors$(\phi,\xi_i,\eta_i,\\g),~i=1,...,s$ satisfying (1.21)-(1.23), then $M$ is called Kenmotsu S-manifold if it satisfies the condition$-$}
\begin{equation}
(\bar{\nabla}_X \phi)Y=\sum_{i=1}^s \{g(\phi X,Y)\xi_i-\eta_i (Y)\phi^2 X\}, \tag{1.26}
\end{equation}
where $\bar{\nabla}$ denotes the Riemannian connection with respect to $g$ on $M$.\\\
\par{For a GKM-manifold $M$,}
\begin{equation}
\bar{\nabla}_X \xi_i =-\phi^2 X. \tag{1.27}
\end{equation}
\par{Cauchy-Riemannian (CR) submanifolds were first introduced in Kaehler geometry, as a bridge between complex and totally real submanifolds. In Bejancu's definition, roughly speaking, their tangent bundle splits into a complex part of constant dimension and a totally real part, orthogonal to the 1st one. after its introduction[1], the definition was soon extended to other ambient spaces and gave rise to a large amount of literature. CR-submanifolds of Kaehler manifolds with Riemannian metric were introduced by A.Bejancu[1][2]. He also introduced totally umbilical CR-submanifolds of Kaehler manifolds. In 1991 and 1992, M.Hasan Shahid extensively studied CR-submanifolds of almost contact manifolds and trans-Sasakian manifolds[6][7][8]. Later in 1993[3], CR-lightlike submanifolds of indefinite Kaehler manifolds was introduced. M.Kobayashi discussed CR-submanifolds of a Sasakian manifold in 1981[11]. In 2010, M.Tarafdar et al. studied contact CR-submanifolds of an indefinite Sasakian manifold. In 2013-2014, B.Laha et al. studied extensively CR-submanifolds of $\epsilon$-paracontact Sasakian manifold, indefinite Lorentzian para-Sasakian manifold[12][13] and submanifolds of some contact and paracontact manifolds[14]. Motivated from their work I have established some new results on CR-submanifolds of $(LCS)_n$-manifold, LP-cosymplectic manifold, S-manifold and Generalized Kenmotsu (GKM) manifold.}\\\
\par{An m-dimensional submanifold $M$ of a differentiable manifold $\tilde{M}$ is called a CR-submanifold if}\\\
$\bullet$ $\xi$ is tangent to $M$,\\\
$\bullet$ tangent bundle splits into two orthogonal differentiable distributions $D$ and $D^\perp$ such that $TM=D\oplus D^\perp$,\\\
$\bullet$ $D:x\longmapsto D_x\subset T_x(M)$, such that $D_x$ is invariant under $\phi$; i.e., $\phi D_x\subset D_x$,\\\
$\bullet$ and for each $x\in M$ and the orthogonal complimentary distribution $D^\perp:x\rightarrow D^\perp_x\subset T_x(M)$ of the distribution $D$ on $M$ is totally real; i.e., $\phi D_x^\perp \subset T^\perp_x(M)$, where $T_x(M)$ and $T^\perp_x(M)$ are the tangent space and normal space respectively at $x$.\\\
\par{We call $D$(respectively $D^\perp$) the horizontal(respectively vertical) distribution. Also CR-submanifold is called $\xi$-horizontal(respectively $\xi$-vertical) if $\xi_x\in D_x$(respectively $\xi_x\in D^\perp$) for each $x\in M$[11].}\\\
\par{The distribution $D$(respectively $D^\perp$) can be defined by a projector $P$ (respectively $Q$), satisfying}
\begin{equation}
P^2=P,~Q^2=Q,~PQ=QP=0,~g\circ(P\times Q)=0. \tag{1.28}
\end{equation}
\par{Let $\bar{M}$ be a differentiable manifold and $M$ be an m-dimensional submanifold, isometrically immersed in $\bar{M}$. Let $g,\nabla$ be the induced metric and Levi-Civita connection on $M$ respectively and $\nabla^\perp,h$ be the normal connection induced by $\bar{\nabla}$ on the normal bundle $TM^\perp$ and the second fundamental form of $M$ respectively.}\\\
\par{Therefore we can decompose the tangent bundle as}
\begin{equation}
T\bar{M}=TM\oplus T^\perp M. \tag{1.29}
\end{equation}
\par{For a vector field $X$ tangent to $M$ and from the fact $\phi D\subset D$, we can write from equation (1.29)[4],}
\begin{equation}
\phi X=PX+QX~ or ~X=PX+QX, \tag{1.30}
\end{equation}
where $PX$(respectively $QX$) is tangential(respectively normal) component of $\phi X$. Authors[10] have defined an endomorphism $P:TM\rightarrow TM$, whose square $P^2$ will be denoted by $Q^\prime$.\\\
\par{Again \begin{equation}
N=BN+CN \tag{1.31}
\end{equation}}
for a vector field $N$ in the normal bundle, where $BN$(respectively $CN$) is a tangential(respectively normal) component of $\phi N$.\\\
\par{Let $TM$ and $T^\perp M$ be the Lie algebras of vector fields tangential and normal to $M$ respectively. $h$ and $A$ denote the 2nd fundamental form and the shape operator of the immersion of $M$ into $\bar{M}$ respectively. If $\nabla$ be the induced connection on $M$, the Gauss and Weingarten formula of $M$ into $\bar{M}$ are characterised by}
\begin{equation}
\bar{\nabla}_X Y=\nabla_X Y+h(X,Y), \tag{1.32}
\end{equation}
\begin{equation}
\bar{\nabla}_X N=-A_N X+\nabla^\perp_X N \tag{1.33}
\end{equation}
for any $X,Y\in TM$ and $N\in T^\perp M$. $\nabla^\perp$ is the connection on normal bundle and $A_N$ is the Weingarten endomorphism associated with $N$ by
\begin{equation}
g(A_N X,Y)=g(h(X,Y),N). \tag{1.34}
\end{equation}
\par{The equation of Gauss is given by}
\begin{equation}
\bar{R}(X,Y,Z,W)=R(X,Y,Z,W)+g(h(X,Z),h(Y,W))-g(h(X,W),h(Y,Z)), \tag{1.35}
\end{equation}
where $\bar{R}$(respectively $R$) is the curvature tensor of $\bar{M}$(respectively $M$).\\\
\par{A submanifold $M$ of a differentiable manifold $\bar{M}$ is called totally umbilical if}
\begin{equation}
h(X,Y)=g(X,Y)H, \tag{1.36}
\end{equation}
where $H$ is the mean curvature denoted by $H=\sum_{i=1}^{k} h(e_i,e_i)$, $k$ is the dimension of $M$ and ${e_1,...,e_k}$ is the local orthonormal frame on $M$.\\\
\par{A submanifold $M$ is called totally geodesic if $h(X,Y)=0$ for each $X,Y\in TM$ and is called minimal if $H=0$ on $M$.}\\\
\par{The covariant derivatives of the tensor fields $\phi,P,Q$ are defined as$-$}
\begin{equation}
(\bar{\nabla}_X \phi)Y= \bar{\nabla}_X\phi Y-\phi\bar{\nabla}_X Y~~\forall X,Y\in T\bar{M}, \tag{1.37}
\end{equation}
and \begin{equation}
(\bar{\nabla}_X P)Y=\nabla_X PY-P(\nabla_X Y), \tag{1.38}
\end{equation}
\begin{equation}
(\bar{\nabla}_X Q)Y=\nabla^\perp_X QY-Q(\nabla_X Y)~~\forall X,Y\in TM. \tag{1.39}
\end{equation}
\par{A CR-submanifold $M$ of a differentiable manifold $\bar{M}$ is called $D$-totally geodesic ($D^\perp$-totally geodesic) if[11]
\begin{equation}
h(X,Y)=0~~~\forall X,Y\in D~~(\forall X,Y\in D^\perp). \tag{1.40}
\end{equation}
\par{A CR-submanifold $M$ of a differentiable manifold $\bar{M}$ is called $D$-umbilic ($D^\perp$-umbilic) if $\forall X,Y\in D$ ($\forall X,Y\in D^\perp$),$L$ being some normal vector field[11],
\begin{equation}
h(X,Y)=g(X,Y)L. \tag{1.41}
\end{equation}
\par{A contact CR-submanifold $M$ of a differentiable manifold $\bar{M}$ is called mixed totally geodesic if[11]
\begin{equation}
h(X,Y)=0~~\forall X\in D~and~\forall Y\in D^\perp. \tag{1.42}
\end{equation}
\par{A contact CR-submanifold $M$ of a differentiable manifold $\bar{M}$ is called foliate contact CR-submanifold if $D$ is involute. $\qquad$ $\qquad \qquad \qquad$  (1.43)}\\\
\par{If $M$ is a foliate $\xi$-horizontal contact CR-submanifold, then the second fundamental form satisfies$-$}
\begin{equation}
h(\phi X,\phi Y)=h(\phi ^2 X,Y)=-h(X,Y). \tag{1.44}
\end{equation}
\par{A submanifold $M$ of an almost contact manifold $\bar{M}$ is called an contact CR-product[9] if it is locally a Riemannian product of $M^T$ and $M^\perp$, where $M^T$ and $M^\perp$ are the leaves of the distributions $D$ and $D^\perp$ respectively.}\\\
$~~~~~~~~~~~~~~~~~~~~~~~~~~~~~~~~~$\section{\large CR-submanifolds of $(LCS)_n$-manifold, LP-cosymplectic manifold, S-manifold and Generalized Kenmotsu (GKM) manifold }
This chapter consists of four sections which are devoted to the study of CR-submanifolds of $(LCS)_n$-manifold, LP-cosymplectic manifold, S-manifold and Generalized Kenmotsu (GKM) manifold.\\\
\subsection{\large CR-submanifolds of an $(LCS)_n$-manifold}
In this section I have obtained some results on CR-submanifolds of an $(LCS)_n$-manifold and also I have given an example of a CR-submanifold of an $(LCS)_n$-manifold. They are as follows$-$\\\
\\\
\\\
\textbf{Theorem 2.1.1.} \textit{Let $M$ be a totally umbilical contact CR-submanifold of an $(LCS)_n$-manifold $\tilde{M}$. Then the anti-invariant distribution $D^\perp$ is 1-dimensional i.e., dim$D^\perp$=1.}\\\
\\\
\textit{Proof.} For an $(LCS)_n$-manifold, from (1.16) we have $\forall X,Y\in D^\perp$,\\\
$(\bar{\nabla}_X \phi)Y=\alpha [g(X,Y)\xi + 2\eta(X)\eta(Y)\xi + \eta(Y)X]$\\\
$\Rightarrow \bar{\nabla}_X \phi Y= \phi \bar{\nabla}_X Y+\alpha[g(X,Y)\xi + 2\eta(X)\eta(Y)\xi + \eta(Y)X]$.\\\
\par{Since $M$ is totally umbilical, by Weingarten and Gauss formula we have,}\\\
$\nabla^\perp_X \phi Y -g(H,\phi Y)X = \phi [\nabla_X Y + g(X,Y)H] + \alpha [g(X,Y)\xi + 2\eta(X)\eta(Y)\xi\\ + \eta(Y)X]$.\\\
\par{Taking inner product with $X$ we get,}\\\
$-g(H,\phi Y)\|X\|^2= g(X,Y)g(\phi H,X)+\alpha g(X,Y)\eta(X)+\alpha2\eta(X)\eta(Y)\eta(X)\\\
+ \alpha \eta(Y)\|X\|^2.\qquad \qquad \qquad \qquad \qquad \qquad \qquad \qquad \qquad \qquad \qquad \qquad$(2.1)
\par{Interchanging $X,Y$ above we get,}\\\
$-g(H,\phi X)\|Y\|^2= g(X,Y)g(\phi H,Y)+\alpha g(X,Y)\eta(Y)+\alpha2\eta(Y)\eta(X)\eta(Y)+\\ \alpha \eta(X)\|Y\|^2.\qquad \qquad \qquad \qquad \qquad \qquad \qquad \qquad \qquad \qquad \qquad \qquad \quad$(2.2)\\\
\par{Using (2.2) in (2.1) and after some calculations we get,}\\\
$g(H,\phi Y)\Big[1-\frac{g(X,Y)^2}{\|X\|^2\|Y\|^2}\Big]+\alpha \eta(Y)\Big[1-\frac{g(X,Y)^2}{\|X\|^2\|Y\|^2}\Big]=\frac{2\alpha}{\|X\|^2}\eta(X)\eta(Y)\Big[\frac{g(X,Y)g(Y,\xi)}{g(Y,Y)}\\-g(X,\xi)\Big]$\\\
...it has a solution if $X\|Y$ i.e., dim$D^\perp$=1.\\\
\\\
\\\
\textbf{Corollary 2.1.1.} \textit{A contact CR-submanifold of an $(LCS)_n$-manifold reduces to a cosymplectic manifold provided the vector field $Z$ becomes the structure vector field $\xi$.}\\\
\\\
\textit{Proof.} In a contact CR-submanifold of an $(LCS)_n$-manifold,\\\
$~~~~~~~~~~~~~~~~$ $(\bar{\nabla}_Z \phi)W=\alpha[g(Z,W)\xi + 2\eta(Z)\eta(W)\xi + \eta(W)Z]$.\\\
\par{Taking inner product with $Z$ we get,}\\\
$g((\bar{\nabla}_Z \phi)W,Z)=\alpha[g(Z,W)\eta(Z)+ 2\eta^2(Z)\eta(W)+ \eta(W)\|Z\|^2]$.\\\
\par{When $Z=\xi$ we have,}\\\
$g((\bar{\nabla}_\xi \phi)W,\xi)=\alpha[-\eta(W)+ 2\eta(W)-\eta(W)]=0$.\\\
\par{Hence the manifold becomes cosymplectic provided $Z=\xi$.}\\\
\\\
\\\
\textbf{Proposition 2.1.1.} \textit{Let $M$ be a CR-submanifold of an $(LCS)_n$-manifold. Then $M$ is $D$-totally geodesic if and only if $A_N X\in D^\perp~~\forall X\in D$ and for any normal vector field $N$.}\\\
\\\
\textit{Proof.} Let $M$ be $D$-totally geodesic. Then $h(X,Y)=0~~\forall X,Y\in D$, so,
$0=g(h(X,Y),N)=g(A_N X,Y)\Rightarrow A_N X\in D^\perp$.\\\
\par{Conversely, let $A_N X\in D^\perp~~\forall X\in D$. Then $\forall Y\in D$,}\\\
$~~~~~~~~~~~~~~~~~~~~~~~~~~~$ $g(A_N X,Y)=g(h(X,Y),N)=0$\\\
$~~~~~~~~~~~~~~~~~~~~~~~~~~~$ $\Rightarrow h(X,Y)=0$ $\forall X,Y\in D$\\\
$~~~~~~~~~~~~~~~~~~~~~~~~~~~$ $\Rightarrow M$ is $D$-totally geodesic.\\\
\\\
\\\
\textbf{Proposition 2.1.2.} \textit{Let $M$ be a CR-submanifold of an $(LCS)_n$-manifold $\tilde{M}$. Then $M$ is $D^\perp$-totally geodesic if and only if $A_N X\in D~~\forall X\in D^\perp$ and for any normal vector field $N$.}\\\
\\\
\textit{Proof.} Similar to Proposition 2.1.1.\\\
\\\
\\\
\textbf{Theorem 2.1.2.} \textit{Let $M$ be a CR-submanifold of an $(LCS)_n$-manifold. If $M$ is $\xi$-horizontal then the distribution $D$ is integrable if and only if}\\\
$~~~~~~~~~~~~~~~~~~~~~~~~~~~$ $h(X,\phi Y)=h(\phi X,Y)~~~\forall X,Y\in D$.\\\
\textit{If $M$ is $\xi$-vertical then the distribution $D^\perp$ is integrable if and only if}\\\
$~~~~~~~~~~~~~~~~~~~$ $A_{\phi X}Y-A_{\phi Y}X=\alpha[\eta(Y)X-\eta(X)Y]~~~\forall X,Y\in D^\perp$.\\\
\\\
\textit{Proof.} $(\bar{\nabla}_X \phi)Y=\alpha [g(X,Y)\xi + 2\eta(X)\eta(Y)\xi + \eta(Y)X]$\\\
$~~~~~~~~~~~~$ $\Rightarrow\bar{\nabla}_X \phi Y-\phi(\bar{\nabla}_X Y)=\alpha [g(X,Y)\xi + 2\eta(X)\eta(Y)\xi + \eta(Y)X]$. $\qquad \qquad$(2.3)\\\
\par{$\bar{\nabla}_X \phi Y$}\\\
$=\bar{\nabla}_X(\phi PY+\phi QY)$\\\
$=\nabla_X \phi PY+ h(X,\phi PY)-A_{\phi QY}X+ \nabla^\perp_X \phi QY$\\\
$=(P\nabla_X \phi PY+ Q\nabla_X \phi PY)- PA_{\phi QY}X- QA_{\phi QY}X+ \nabla^\perp_X \phi QY$\\\
$+ h(X,\phi PY).$ $\quad \qquad \qquad \qquad \qquad \qquad \qquad \qquad \qquad \qquad \qquad \qquad \quad$ (2.4)\\\
\par{$\phi(\bar{\nabla}_X Y)$}\\\
$=\phi(\nabla_X Y+h(X,Y))$\\\
$=\phi(P\nabla_X Y+Q\nabla_X Y)+Bh(X,Y)+ Ch(X,Y)$.$\qquad \qquad \qquad \qquad \qquad$(2.5)\\\
\par{Using (2.4), (2.5) on (2.3) we have,}\\\
$P\nabla_X \phi PY+Q\nabla_Y\phi PY-PA_{\phi QY}X-QA_{\phi QY}X+\nabla^\perp_X \phi QY+h(X,\phi PY)\\-\phi(P\nabla_X Y)-\phi(Q\nabla_X Y)-Bh(X,Y)-Ch(X,Y)= \alpha[g(X,Y)P\xi\\+g(X,Y)Q\xi+2\eta(X)\eta(Y)P\xi+2\eta(X)\eta(Y)Q\xi+\eta(Y)(PX+QX)]$.\\\
\par{Comparing horizontal, vertical and normal parts we obtain,}\\\
$P\nabla_X \phi PY-PA_{\phi QY} X-\phi P\nabla_XY=\alpha g(X,Y)P\xi+\alpha\eta(Y)PX\\+2\alpha\eta(X)\eta(Y)P\xi$,\\\
$Q\nabla_X\phi PY-QA_{\phi QY}X-Bh(X,Y)=\alpha g(X,Y)Q\xi+\alpha\eta(Y)QX\\+2\alpha\eta(X)\eta(Y)Q\xi$,\\\
$h(X,\phi PY)+\nabla^\perp_X \phi QY-\phi(Q\nabla_XY)-Ch(X,Y)=0$.$\qquad \qquad \qquad \qquad$(2.6)\\\
\par{If $M$ is $\xi$-horizontal, from (2.6) we get,}
\begin{equation}
h(X,\phi PY)=\phi(\nabla_X Y)+Ch(X,Y)~~~\forall X,Y\in D. \tag{2.7}
\end{equation}
\par{Interchanging $X,Y$ we have,}
\begin{equation}
h(Y,\phi PX)=\phi(\nabla_Y X)+Ch(X,Y)~~~\forall X,Y\in D. \tag{2.8}
\end{equation}
\par{Subtracting (2.8) from (2.9) we get,}\\\
$h(X,\phi PY)-h(Y,\phi PX)=\phi(Q[X,Y])=0$ if and only if $[X,Y]\in D$ if and only if $D$ is integrable. Hence $D$ is integrable if and only if $h(X,\phi Y)=h(Y,\phi X)~~~\forall X,Y\in D$.\\\
\
\par{If $M$ is $\xi$-vertical, from (2.6) we have, $\forall X,Y\in D^\perp$,}
\begin{equation}
\nabla^\perp_X\phi Y=\phi(Q\nabla_X Y)+Ch(X,Y), \tag{2.9}
\end{equation}
thus, $\nabla^\perp_X\phi Y=A_{\phi Y}X+\bar{\nabla}_X \phi Y=A_{\phi Y}X+(\bar{\nabla}_X \phi)Y+\phi(\bar{\nabla}_X Y)\\=A_{\phi Y}X+\alpha [g(X,Y)\xi + 2\eta(X)\eta(Y)\xi + \eta(Y)X]+\phi(\nabla_X Y+h(X,Y))\\=A_{\phi Y}X+\alpha [g(X,Y)\xi + 2\eta(X)\eta(Y)\xi + \eta(Y)X]+\phi P\nabla_X Y+\phi Q\nabla_X Y\\+Bh(X,Y)+Ch(X,Y).$\\\
\par{Now from (2.9) we get,}\\\
$\phi (Q\nabla_X Y)+Ch(X,Y)=A_{\phi Y}X+\alpha[g(X,Y)\xi+2\eta(X)\eta(Y)\xi+\eta(Y)X]\\+\phi P\nabla_X Y+\phi Q\nabla_X Y+Bh(X,Y)+Ch(X,Y)$\\\
$\Rightarrow\phi P\nabla_X Y=-A_{\phi Y}X-\alpha[g(X,Y)\xi+2\eta(X)\eta(Y)\xi+\eta(Y)X]\\-Bh(X,Y). \qquad \qquad \qquad \qquad \qquad \qquad \qquad \qquad \qquad \qquad \qquad \qquad $ (2.10)\\\
\par{Interchanging $X,Y$ we get,}\\\
$\phi P\nabla_Y X=-A_{\phi X}Y-\alpha[g(X,Y)\xi+2\eta(X)\eta(Y)\xi+\eta(X)Y]-Bh(X,Y)$.$~~~~~~~~~~~~~~~~~~~~~~~~~~~~~~~~~~~~~~~~~~$(2.11)\\\
\par{Subtracting (2.11) from (2.10) we obtain,}\\\
$\phi P([X,Y])=-A_{\phi Y}X+A_{\phi X}Y-\alpha[\eta(Y)X-\eta(X)Y]$.\\\
\par{Hence $[X,Y]\in D^\perp$ if and only if $A_{\phi X}Y-A_{\phi Y}X=\alpha[\eta(Y)X-\eta(X)Y]\\\forall X,Y\in D^\perp$.}\\\
\\\
\\\
\textbf{Proposition 2.1.3.} \textit{Let $M$ be a $D$-umbilic CR-submanifold of an $(LCS)_n$-manifold $\tilde{M}$. If $M$ is $\xi$-horizontal(respectively $\xi$-vertical), then $M$ is $D$-totally geodesic(respectively $D^\perp$-totally geodesic) CR-submanifold.}\\\
\\\
\textit{Proof.} Consider $M$ as $D$-umbilic $\xi$-horizontal CR-submanifold. Then from definition 1.41 we have, $h(X,Y)=g(X,Y)L~~~\forall X,Y\in D$, L being some normal vector field on $M$.\\\
\par{By putting $X=Y=\xi$ and as $h(\xi,\xi)=0$ we have,}\\\
$g(\xi,\xi)L=h(\xi,\xi)=0$\\\
$\Rightarrow L=0$ (since $g(\xi,\xi)=\eta(\xi)=-1$),\\\
hence $h(X,Y)=g(X,Y)L=0$\\\
$\Rightarrow M$ is a $D$-totally geodesic CR-submanifold.\\\
\par{Similarly, it can be easily shown that if $M$ is $D^\perp$-umbilic $\xi$-vertical CR-submanifold, then it is $D^\perp$-totally geodesic CR-submanifold.}\\\
\\\
\\\
\textbf{Lemma 2.1.1.} \textit{Let $M$ be a CR-submanifold of an $(LCS)_n$-manifold. Then $M$ is mixed totally geodesic if and only if $A_N X\in D~~\forall X\in D$ and $A_NX\in D^\perp~~\forall X\in D^\perp$, $N$ is any normal vector field.}\\\
\\\
\textit{Proof.} If $M$ is mixed totally geodesic, then from (1.42) we have $\forall X\in D,\\Y\in D^\perp$,\\\
$h(X,Y)=0$\\\
$\Rightarrow g(h(X,Y),N)=0$\\\
$\Rightarrow g(A_N X,Y)=0~~\forall Y\in D^\perp$ (by (1.34))\\\
$\Rightarrow A_N X\in D~~\forall X\in D$.\\\
\par{Similarly, $A_NX\in D^\perp~~\forall X\in D^\perp$.}\\\
\par{Converse can be proved easily.}
\\\
\\\
\textbf{Theorem 2.1.3.} \textit{If $M$ is a mixed totally geodesic CR-submanifold of an $(LCS)_n$-manifold, then $\forall X\in D$\textit{ and for any normal vector field} $N$,}\\\
$A_{\phi N}X=\phi A_NX$ and $\nabla^\perp_X \phi N=\phi\nabla^\perp_XN$.\\\
\\\
\textit{Proof.} From Lemma 2.1.1 and using (1.32) we have,\\\
$~~~~~~~~~~~~~~~~~~~~~~~~$ $\nabla_X\phi N=\bar{\nabla}_X\phi N$ (since $h(X,\phi N)=0$)\\\
$~~~~~~~~~~~~~~~~~~~~~~~~~~~~~~~~~~$ $=(\bar{\nabla}_X\phi)N+\phi(\bar{\nabla}_X N)$\\\
$~~~~~~~~~~~~~~~~~~~~~~~~~~~~~~~~~~$ $=(\bar{\nabla}_X\phi)N+\phi(\nabla^\perp_X N-A_N X)$\\\
$~~~~~~~~~~~~~~~~~~~~~~~~~~~~~~~~~~$ $=\phi(\nabla^\perp_X N-A_N X)$.\\\
\par{Again $\nabla_X \phi N=-A_{\phi N}X+\nabla^\perp_X \phi N$.}\\\\
\par{Comparing above two equations we get $A_{\phi N}X=\phi A_NX$ and $\nabla^\perp_X \phi N=\phi\nabla^\perp_XN$.}\\\
\\\
\\\
\textbf{Proposition 2.1.4.} \textit{If $M$ is a foliate $\xi$-horizontal mixed totally geodesic CR-submanifold of an $(LCS)_n$-manifold, then $\phi A_N X=A_N \phi X~~\forall X\in D$ and for any normal vector field $N$.}\\\
\\\
\textit{Proof.} Using (1.34) and (1.12) we have,
\begin{equation}
g(h(X,\phi Y),N)=g(A_N X,\phi Y)=g(\phi A_N X,Y). \tag{2.12}
\end{equation}
\par{Also using (1.34) we have, $g(h(\phi X,Y),N)=g(A_N \phi X,Y).\qquad \quad$ (2.13)}\\\
\par{Since $M$ is $\xi$-horizontal foliate CR-submanifold we have, $h(X,\phi Y)=h(\phi X,Y)$ and so from (2.12) and (2.13) we get, $\phi A_N X=A_N \phi X$.}\\\
\\\
\\\
\textbf{Theorem 2.1.4.} \textit{Let $M$ be a $\xi$-horizontal CR-submanifold of an $(LCS)_n$-manifold $\tilde{M}$. Then the leaf $M^T$ of $D$ is totally geodesic in $M$ if and only if $g(h(D,D),\phi D^\perp)=0$.}\\\
\\\
\textit{Proof.} Let the leaf $M^T$ of $D$ is totally geodesic in $M$, then for $X,Y\in D$, $\nabla_X \phi Y\in D$.\\\
\par{Now $(\bar{\nabla}_X \phi)Y=\bar{\nabla}_X \phi Y-\phi(\bar{\nabla}_X Y)$}\\\
$=\bar{\nabla}_X PY+\bar\nabla_X QY-\phi(\nabla_X Y+h(X,Y))$\\\
$~~~~~~~~~~~~~~~~~~~~~~~~$ $=\nabla_X PY+h(X,PY)-A_{QY}X+\nabla_X^\perp QY-P(\nabla_X Y)-Q(\nabla_X Y)-Bh(X,Y)$\\\
$-Ch(X,Y)$\\\
$~~~~~~~~~~~~~~~~~~~~~~~~$ $=(\bar{\nabla}_X P)Y+(\bar{\nabla}_X Q)Y+h(X,PY)-A_{QY}X-Bh(X,Y)-Ch(X,Y)$.$~~~$(2.14)\\\
\par{Now from (1.16) we have,}\\\
$(\bar\nabla_X P)Y+(\bar{\nabla}_X Q)Y+h(X,PY)-A_{QY}X-Bh(X,Y)-Ch(X,Y)\\=\alpha[g(X,Y)\xi+2\eta(X)\eta(Y)\xi+\eta(Y)X]$\\\
$\Rightarrow(\bar{\nabla}_X P)Y=A_{QY}X+Bh(X,Y)+\alpha g(X,Y)\xi+\alpha\eta(Y)PX+2\alpha\eta(X)\eta(Y)\xi$,\\\
and $(\bar{\nabla}_X Q)Y=Ch(X,Y)-h(X,PY)+\alpha\eta(Y)QX$\\\
$\Rightarrow P\nabla_X Y=\nabla_XPY-A_{QY}X-Bh(X,Y)-\alpha g(X,Y)\xi-\alpha\eta(Y)PX\\-2\alpha\eta(X)\eta(Y)\xi$, $\qquad \qquad \qquad \qquad \qquad \qquad \qquad \qquad \qquad \qquad \qquad$ (2.15)\\\
and $Q\nabla_X Y=\nabla^\perp_X QY+h(X,PY)-\alpha\eta(Y)QX-Ch(X,Y)$.$\qquad \qquad$(2.16)\\\
\par{For $Z\in D^\perp$, since $\nabla_X \phi Y\in D$,}\\\
$0=g(\nabla_X \phi Y,Z)=-g(\nabla_X Z,\phi Y)=-g(\phi\nabla_X Z,Y)=-g(P\nabla_X Z,Y)$\\\
$=g(A_{QZ}X,Y)+g(Bh(X,Z),Y)$(since $P\nabla_X Z=-A_{QZ}X-Bh(X,Z)$)\\\
$=g(h(X,Y),QZ)$ (using (1.34))\\\
$=g(h(X,Y),\phi Z)$.\\\
\par{Hence $g(h(D,D^\perp)$,$\phi D^\perp)=0$.}\\\
\\\
\par{Conversely, $\forall X,Y\in D, Z\in D^\perp$,$0=g(h(X,\phi Y),\phi Z)=g(\bar{\nabla}_X\phi Y, \phi Z)=g(\phi \bar{\nabla}_XY,\phi Z)=g(\bar{\nabla}_XY,Z)=g(\nabla_XY,Z).$}\\\
\par{Hence $\nabla_XY\in D~~\forall X,Y\in D$ which means that the leaf $M^T$ of $D$ is totally geodesic in $M$.}
\\\
\\\
\textbf{Theorem 2.1.5.} \textit{Let $M$ be a CR-submanifold of an $(LCS)_n$-manifold $\tilde{M}$. Then the leaf $M^\perp$ of $D^\perp$ is totally geodesic in $M$ if and only if}\\\
$~~~~~~~~~~~~~~~~~$ $g(h(W,Y),\phi Z)+\alpha\eta(Y)g(W,Z)+2\alpha\eta(Y)\eta(Z)\eta(W)=0~~\forall Y\in D$ and $\forall W,Z\in D^\perp$.\\\
\\\
\textit{Proof.} From (2.15) we have on putting $X=W,Y=Z\in D^\perp$,\\\
\centerline{$P\nabla_W Z=-A_{QZ}W-Bh(W,Z)-\alpha g(W,Z)\xi-2\alpha\eta(W)\eta(Z)\xi$.}\\\
\par{Taking inner product with $Y\in D$ we get,}\\\
$g(P\nabla_W Z,Y)=-g(A_{QZ}W,Y)-g(Bh(W,Z),Y)-\alpha g(W,Z)\eta(Y)\\-2\alpha\eta(W)\eta(Z)\eta(Y)$.\\\
\par{Now, $g(P\nabla_W Z,Y)=g(\phi\nabla_WZ,Y)=g(\nabla_W Z,\phi Y)=g(\nabla_W Z,PY)$.\\\
\par{Thus, $g(\nabla_W Z,PY)=-g(h(W,Y),QZ)-\alpha\eta(Y)g(W,Z)\\-2\alpha\eta(W)\eta(Z)\eta(Y)=-[g(h(W,Y),\phi Z)+\alpha\eta(Y)g(W,Z)\\+2\alpha\eta(Y)\eta(Z)\eta(W)].$}\\\
\par{Hence the leaf $M^\perp$ of $D^\perp$ is totally geodesic in $M$ if and only if $\forall Y\in D$ and $\forall W,Z\in D^\perp,g(h(W,Y),\phi Z)+\alpha\eta(Y)g(W,Z)+2\alpha\eta(Y)\eta(Z)\eta(W)=0$.\\\
\\\
\\\
\textbf{Corollary 2.1.2.} \textit{Let $M$ be a $\xi$-vertical CR-submanifold of an $(LCS)_n$-manifold $\tilde{M}$. Then the leaf $M^\perp$ of $D^\perp$ is totally geodesic in $M$ if and only if $g(h(D^\perp,D),\phi D^\perp)=0$.}\\\
\\\
\\\
\textbf{Theorem 2.1.6.} \textit{Let $M$ be a $\xi$-horizontal CR-submanifold of an $(LCS)_n$-manifold $\tilde{M}$. Then $M$ is a contact CR-product if and only if $A_{\phi W}X+\alpha\eta(X)W=0~~\forall X\in D,W\in D^\perp$.}\\\
\\\
\textit{Proof.} Let $M$ be a contact CR-product. From Theorem 2.1.5 we have,\\\
$~~~~~~~~~~~~~~~~~~$ $g(A_{\phi Z}W,Y)+\alpha\eta(Y)g(W,Z)=0,\forall Z,W\in D^\perp, X\in D$ (since $M$ is $\xi$-horizontal, $\eta(Z)=0,\eta(W)=0$)\\\
$~~~~~~~~~~~~~~~~~~$ $\Rightarrow g(h(W,Y),\phi Z)+\alpha\eta(Y)g(W,Z)=0$\\\
$~~~~~~~~~~~~~~~~~~$ $\Rightarrow g(A_{\phi Z}Y,W)+\alpha\eta(Y)g(W,Z)=0$\\\
$~~~~~~~~~~~~~~~~~~$ $\Rightarrow g(A_{\phi W}X+\alpha\eta(X)W,Z)=0$ (interchanging $Z,W$ and replacing $Y$ by $X$)\\\
$~~~~~~~~~~~~~~~~~~$ $\Rightarrow A_{\phi W}X+\alpha\eta(X)W\in D$.$\qquad \qquad \qquad \qquad \qquad \qquad \qquad \qquad \qquad$(2.17)\\\
\par{As the leaf $M^T$ of $D$ is totally geodesic in $M$, for $Y\in D$,}\\\
$~~~~~~~~~~~~~~~~~~~~~~$ $g(A_{\phi W}X+\alpha\eta(X)W,Y)$\\\
$~~~~~~~~~~~~~~~~~~~~~~$ $=g(A_{\phi W}X,Y)+\alpha\eta(X)g(W,Y)$\\\
$~~~~~~~~~~~~~~~~~~~~~~$ $=g(h(X,Y),\phi W)$\\\
$~~~~~~~~~~~~~~~~~~~~~~$ $=g(\phi h(X,Y),W)$\\\
$~~~~~~~~~~~~~~~~~~~~~~$ $=g(\phi\bar{\nabla}_X Y,W)$\\\
$~~~~~~~~~~~~~~~~~~~~~~$ $=g(\bar{\nabla}_X\phi Y,W)$\\\
$~~~~~~~~~~~~~~~~~~~~~~$ $=g(\nabla_X\phi Y,W)$\\\
$~~~~~~~~~~~~~~~~~~~~~~$ $=0$\\\
$~~~~~~~~~~~~~~~~~~~~~~$ $\Rightarrow A_{\phi W} X+\alpha\eta(X)W\in D^\perp$.$\qquad \qquad \qquad \qquad \qquad \qquad \qquad \qquad \qquad$(2.18)\\\
\par{From (2.17) and (2.18) we get, $A_{\phi W}X+\alpha\eta(X)W=0$.}\\\
\par{Conversely, $\forall X\in D,W,Z\in D^\perp$,}\\\
$~~~~~~~~~~~~~~~~~~~$ $g(A_{\phi W}X,Z)+\alpha\eta(X)g(W,Z)=0$\\\
$~~~~~~~~~~~~~~~~~~~~~~$ $\Rightarrow g(h(X,Z),\phi W)+\alpha\eta(X)g(W,Z)+2\alpha\eta(X)\eta(Z)\eta(W)=0$ (since $M$ is $\xi$-horizontal, $\eta(Z)=0,\eta(W)=0$)\\\
$~~~~~~~~~~~~~~~~~~$ $\Rightarrow$ leaf $M^\perp$ of $D^\perp$ is totally geodesic in $M$ (by Theorem 2.1.5).\\\
\\\
$~~~~~~~~~~~~~~~~~~~~$ $g(\nabla_X Y,Z)=g(\bar{\nabla}_X Y,Z)$ (by (1.32))\\\
$~~~~~~~~~~~~~~~~~~~~~~~~~~~~~~~~~~~$ $=g(\phi \bar{\nabla}_X Y,\phi Z)$\\\
$~~~~~~~~~~~~~~~~~~~~~~~~~~~~~~~~~~~$ $=g(\bar{\nabla}_X \phi Y,\phi Z)$\\\
$~~~~~~~~~~~~~~~~~~~~~~~~~~~~~~~~~~~$ $=g(h(X,\phi Y),\phi Z)$ (by (1.32))\\\
$~~~~~~~~~~~~~~~~~~~~~~~~~~~~~~~~~~~$ $=g(A_{\phi Z}X,\phi Y)$ (by (1.34))\\\
$~~~~~~~~~~~~~~~~~~~~~~~~~~~~~~~~~~~$ $=-\alpha\eta(X)g(\phi Y,Z)$ (by the condition)\\\
$~~~~~~~~~~~~~~~~~~~~~~~~~~~~~~~~~~~$ $=0$\\\
$~~~~~~~~~~~~~~~~~~~~$ $\Rightarrow\nabla_XY\in D~~\forall X,Y\in D$\\\
$~~~~~~~~~~~~~~~~~~~~$ $\Rightarrow$ leaf $M^T$ of $D$ is totally geodesic in $M$.\\\
\par{Hence $M$ is a contact CR-product.}\\\
\\\
\\\
\textbf{Theorem 2.1.7.} \textit{Let $M$ be a $M$ be a $\xi$-horizontal contact CR-product of an $(LCS)_n$-manifold $\tilde{M}$. Then for unit vectors $X\in D$ and $Z\in D^\perp$ with $\eta(X)=0$,}\\\
i)$g(h(\nabla_X \phi X,Z),\phi Z)=\alpha^2,$\\\
ii)$g(h(\nabla_{\phi X}X,Z),\phi Z)=\alpha^2,$\\\
iii)$g(h(\phi X,\nabla_X Z),\phi Z)=0,$\\\
iv)$g(h(X,\nabla_{\phi X}Z),\phi Z)=0$.\\\
\\\
\textit{Proof.} i)As $M$ is a $\xi$-horizontal contact CR-product of an $(LCS)_n$-manifold $\tilde{M}$, from Theorem 2.1.6 we have,\\\
$g(h(\nabla_X \phi X,Z),\phi Z)=-\alpha\eta(\nabla_X \phi X)g(Z,Z)$\\\
$~~~~~~~~~~~~~~~~~~~~~~~~~~~$ $=-\alpha g(\nabla_X \phi X,\xi)$\\\
$~~~~~~~~~~~~~~~~~~~~~~~~~~~$ $=\alpha g(\phi X,\nabla_X \xi)$\\\
$~~~~~~~~~~~~~~~~~~~~~~~~~~~$ $=\alpha g(\phi X,\alpha\phi X)$\\\
$~~~~~~~~~~~~~~~~~~~~~~~~~~~$ $=\alpha^2g(\phi X,\phi X)$\\\
$~~~~~~~~~~~~~~~~~~~~~~~~~~~$ $=\alpha^2g(X,X)$\\\
$~~~~~~~~~~~~~~~~~~~~~~~~~~~$ $=\alpha^2$.\\\
\par{In the same way we can obtain the rest three equations.}\\\
\\\
\\\
\\\
$\bullet$ \underline{\textbf{Example of a CR-submanifold of an $(LCS)_n$-manifold}} \textbf{:}\\\
\
\par{Consider $\bar{M}=\mathbb{R}^7$ be the semi-Euclidean space. Define the semi-Euclidean metric\\ $g=[-dt^2+\sum^6_{i=1}(dx_i)^2],(t,x_1,x_2,...,x_6)\in \bar{M}$. Take $\xi=\frac{\partial}{\partial t},\eta=dt,\phi\big(\frac{\partial}{\partial t}\big)=0$, and}\\\
$\phi\big(\frac{\partial}{\partial x_1}\big)=-\frac{\partial}{\partial x_4},\phi\big(\frac{\partial}{\partial x_2}\big)=-\frac{\partial}{\partial x_5},\phi\big(\frac{\partial}{\partial x_3}\big)=-\frac{\partial}{\partial x_6},\phi\big(\frac{\partial}{\partial x_4}\big)=-\frac{\partial}{\partial x_1},\phi\big(\frac{\partial}{\partial x_5}\big)=-\frac{\partial}{\partial x_2},\\\phi\big(\frac{\partial}{\partial x_6}\big)=-\frac{\partial}{\partial x_3}.$\\\
Then we can see that $(\phi,\xi,\eta,g)$ is an $(LCS)_7$ structure on $\bar{M}$ and hence $\bar{M}$ is an $(LCS)_7$-manifold. Define a submanifold $M$ of $\bar{M}$ by $M=\{(0,x_2,0,x_4,\\x_5,x_6,t)\in\mathbb{R}^7\}$ endowed with the global vector fields\\\
\par{$h_1=\frac{\partial}{\partial x_4},h_2=\xi=\frac{\partial}{\partial t},h_3=\frac{\partial}{\partial x_1}+x_4\frac{\partial}{\partial t},h_4=\frac{\partial}{\partial x_6},h_5=\frac{\partial}{\partial x_2}+x_6\frac{\partial}{\partial t}$.}\\\
Then the distributions $D_T=span\{h_1,h_3\}$ and $D^\perp=span\{h_4,h_5\}$ are respectively invariant and anti-invariant distributions on $\bar{M}$. Thus we can write $TM=D_T(=\{h_1,h_3\})\oplus D^\perp(=\{h_4,h_5\})\oplus\\<h_2(=\{\xi\})>.$ Consequently $M$ is a CR-submanifold of $\bar{M}=\mathbb{R}^7.$\\\
\
\subsection{\large CR-submanifolds of an LP-cosymplectic manifold}
Here I have considered a CR-submanifold of an LP-cosymplectic manifold and using it's structure equation I have stated and proved the following results$-$\\\
\\\
\\\
\textbf{Theorem 2.2.1.} \textit{Let $M$ be a totally umbilical contact CR-submanifold of an LP-cosymplectic manifold $\tilde{M}$. Then the anti-invariant distribution $D^\perp$ is 1-dimensional, i.e., dim$D^\perp$=1.}\\\
\\\
\textit{Proof.} From (1.19) we have $\forall Z,W\in D^\perp$,\\\
$(\bar{\nabla}_Z\phi)W=g(Z,W)\xi+2\eta(Z)\eta(W)\xi=0$\\\
$\Rightarrow\bar{\nabla}_Z\phi W=\phi(\bar{\nabla}_Z W)$.\\\
\par{Since $M$ is totally umbilical, by Weingarten and Gauss formula we have,}\\\
$\nabla^\perp_Z \phi W-g(H,\phi W)Z=\phi[\nabla_Z W+g(Z,W)H]$.\\\
\par{Taking inner product with $Z$ we get,}
\begin{equation}
-g(H,\phi W)\|Z\|^2=g(Z,W)g(\phi H,Z). \tag{2.19}
\end{equation}
\par{Interchanging $Z,W$ in (2.19) and then using (2.19) on it we get after some calculations,}\\\
\centerline{$g(H,\phi Z)\Big[1-\frac{g(Z,W)^2}{\|Z\|^2\|W\|^2}\Big]=0$}\\\
...it has a solution if $Z||W$ i.e., dim$D^\perp$=1.\\\
\\\
\\\
\textbf{Proposition 2.2.1.} \textit{Let $M$ b e a CR-submanifold of an LP-cosymplectic manifold. Then $M$ is a $D$-totally geodesic CR-submanifold if and only if $A_N X\in D^\perp~~\forall X\in D$ and $N$ is any normal vector field.}\\\
\\\
\textit{Proof.} Similar to Proposition 2.1.1.\\\
\\\
\\\
\textbf{Proposition 2.2.2.} \textit{Let $M$ be a CR-sumanifold of an LP-cosymplectic manifold $\tilde{M}$. Then $M$ is $D^\perp$-totally geodesic if and only if $A_N X\in D\\\forall X\in D^\perp$ and $N$ is any normal vector field.}\\\
\\\
\textit{Proof.} Similar to Proposition 2.1.2.\\\
\\\
\\\
\textbf{Theorem 2.2.2.} \textit{Let $M$ be a CR-submanifold of an LP-cosymplectic manifold. If $M$ is $\xi$-horizontal, then the distribution $D$ is integrable if and only if $h(X,\phi Y)=h(Y,\phi X)$ $\forall X,Y\in D$.}\\\
\textit{If $M$ is $\xi$-vertical, then the distribution $D^\perp$ is integrable if and only if}\\\
$~~~~~~~~~~~~~~~~~~$ $A_{\phi Y}X=A_{\phi X}Y$ $\forall X,Y\in D^\perp$.\\\
\\\
\textit{Proof.} $(\bar{\nabla}_X\phi)Y=0$ (from (1.19))\\\
$~~~~~~~~$ $\Rightarrow\bar{\nabla}_X\phi Y=\phi(\bar{\nabla}_X Y)$\\\
$~~~~~~~~~$ $\Rightarrow P\nabla_X \phi PY+Q\nabla_X \phi QY-PA_{\phi QY}X-QA_{\phi QY}X+\nabla^\perp_X\phi QY+h(X,\phi PY)=\phi(P\nabla_XY+Q\nabla_XY)+Bh(X,Y)+Ch(X,Y)$ (using (2.4), (2.5)).\\\
\par{Comparing horizontal, vertical, normal parts we obtain,}
\begin{equation}
P\nabla_X\phi PY-PA_{\phi QY}X=\phi(P\nabla_XY), \tag{2.20}
\end{equation}
\begin{equation}
Q\nabla_X\phi PY-QA_{\phi QY}X=Bh(X,Y), \tag{2.21}
\end{equation}
\begin{equation}
h(X,\phi PY)+\nabla^\perp_X\phi QY=\phi(Q\nabla_XY)+Ch(X,Y). \tag{2.22}
\end{equation}
\par{If $M$ is $\xi$-horizontal, from (2.22) we get for $X,Y\in D$,}
\begin{equation}
h(X,\phi PY)=\phi(Q\nabla_XY)+Ch(X,Y). \tag{2.23}
\end{equation}
\par{Interchanging $X,Y$ above we have,}
\begin{equation}
h(Y,\phi PX)=\phi(Q\nabla_YX)+Ch(X,Y). \tag{2.24}
\end{equation}
\par{Subtracting (2.24) from (2.23) we get,}\\\
$h(X,\phi PY)-h(Y,\phi PX)=\phi(Q[X,Y])=0$ if and only if $[X,Y]\in D$ if and only if $D$ is integrable.\\\
\par{Hence $D$ is integrable if and only if $h(X,\phi Y)=h(Y,\phi X)$.}\\\
\par{If $M$ is $\xi$-vertical then from (2.22) we have for $X,Y\in D^\perp$,}\\\
$\nabla^\perp_X\phi Y=Ch(X,Y)+\phi Q\nabla_XY.$\\\
$\bar{\nabla}_X\phi Y=(\bar{\nabla}_X\phi)Y+\phi(\nabla_XY)$\\\
$~~~~~~~~~~$ $=0+\phi(P\nabla_XY+Q\nabla_XY)+Bh(X,Y)+Ch(X,Y)$ (using (1.19))\\\
$\Rightarrow\nabla^\perp_X\phi Y=A_{\phi Y}X+\phi(P\nabla_XY)+\phi(Q\nabla_XY)+Bh(X,Y)+Ch(X,Y)$\\\
(using (1.33))\\\
$\Rightarrow Ch(X,Y)+\phi Q\nabla_XY=A_{\phi Y}X+\phi(P\nabla_XY)+\phi(Q\nabla_XY)+Bh(X,Y)\\+Ch(X,Y)$\\\
$\Rightarrow\phi(P\nabla_XY)=-A_{\phi Y}X-Bh(X,Y)$.$\qquad \qquad \qquad \qquad \qquad \qquad \qquad$(2.25)\\\
\par{Interchanging $X,Y$ we have,}
\begin{equation}
\phi(P\nabla_YX)=-A_{\phi X}Y-Bh(X,Y). \tag{2.26}
\end{equation}
\par{Subtracting (2.26) from (2.25) we get,}\\\
\centerline{$\phi P[X,Y]=-A_{\phi Y}X+A_{\phi X}Y.$}\\\
\par{Now, $[X,Y]\in D^\perp$ if and only if $A_{\phi X}Y=A_{\phi Y}X$.}\\\
\\\
\\\
\textbf{Proposition 2.2.3.} \textit{Let $M$ be a $D$-umbilic CR-submanifold of an LP-cosymplectic manifold $\tilde{M}$. If $M$ is $\xi$-horizontal(respectively $\xi$-vertical), then $M$ is $D$-totally geodesic(respectively $D^\perp$-totally geodesic).}\\\
\\\
\textit{Proof.} Same as Proposition 2.1.3.\\\
\\\
\\\
\textbf{Lemma 2.2.1.} \textit{Let $M$ CR-submanifold of an LP-cosymplectic manifold. Then $M$ is mixed totally geodesic if and only if $A_NX\in D~~\forall X\in D$ and $A_NX\in D^\perp~~\forall X\in D^\perp$ and $N$ is any normal vector field.}\\\
\\\
\\\
\textbf{Theorem 2.2.3.} \textit{If $M$ is a mixed totally geodesic CR-submanifold of an LP-cosymplectic manifold, then $\forall X\in D$ and for any normal vector field $N$, $A_{\phi N}X=\phi A_NX$ and $\nabla^\perp_X\phi N=\phi\nabla^\perp_XN$.}\\\
\\\
\textit{Proof.} Same as Theorem 2.1.3.\\\
\\\
\\\
\textbf{Proposition 2.2.4.} \textit{If $M$ be a foliate $\xi$-horizontal mixed totally geodesic CR-submanifold of an LP-cosymplectic manifold, then $\phi A_NX=A_N\phi X\\\forall X\in D$ and for any normal vector field $N$.}\\\
\\\
\textit{Proof.} Same as Proposition 2.2.4.\\\
\\\
\\\
\textbf{Theorem 2.2.4.} \textit{Let $M$ be a $\xi$-horizontal CR-submanifold of an LP-cosymplectic manifold $\tilde{M}$. Then the leaf $M^T$ of $D$ is totally geodesic in $M$ if and only if $g(h(D,D),\phi D^\perp)=0$.}\\\
\\\
\textit{Proof.} Let the leaf $M^T$ of $D$ be totally geodesic in $M$.\\\
\par{Using (2.14) we have,}\\\
$0=(\bar{\nabla}_X\phi)Y=(\bar{\nabla}_XP)Y+(\bar{\nabla}_XQ)Y+h(X,PY)-A_{QY}X-Bh(X,Y)\\-Ch(X,Y)$\\\
$\Rightarrow(\bar{\nabla}_XP)Y=A_{QY}X+Bh(X,Y)$,\\\
and $(\bar{\nabla}_XQ)Y=Ch(X,Y)-h(X,PY)$.\\\
\par{Hence $P\nabla_XY=\nabla_XPY-A_{QY}X-Bh(X,Y)$,$\qquad \qquad \qquad \qquad$(2.27)}\\\
and $Q\nabla_XY=\nabla^\perp_XQY-Ch(X,Y)+h(X,PY)$. $\qquad \qquad \qquad \qquad$ (2.28)\\\
\par{For $X,Y\in D,Z\in D^\perp$, since $\nabla_X\phi Y\in D$,}\\\
$0=g(\nabla_X\phi Y,Z)=-g(\nabla_XZ,\phi Y)=-g(\phi\nabla_XZ,Y)$\\\
$=-g(P\nabla_XZ,Y)=g(A_{QZ}X,Y)+g(Bh(X,Z),Y)$\\\
$=g(h(X,Y),QZ)$ (using (1.34))\\\
$=g(h(X,Y),\phi Z)$.\\\
\par{Hence $g(h(D,D),\phi D^\perp)=0$.}\\\
\par{Converse is same as Theorem 2.1.4.}\\\
\\\
\\\
\textbf{Theorem 2.2.5.} \textit{Let $M$ be a CR-submanifold of an LP-cosymplectic manifold $\tilde{M}$. Then the leaf $M^\perp$of $D^\perp$ is totally geodesic in $M$ if and only if $g(h(D^\perp,D),\phi D^\perp)=0.$}\\\
\\\
\textit{Proof.} From (2.27) we have,\\\
\centerline{$P\nabla_WZ=-A_{QZ}W-Bh(W,Z)~~\forall W,Z\in D^\perp$.}\\\
\par{Taking inner product with $Y\in D$ and using (1.34) we get,}\\\
$g(\nabla_WZ,PY)=-g(h(W,Y),QZ)=-g(h(W,Y),\phi Z)\Rightarrow$the leaf $M^\perp$of $D^\perp$ is totally geodesic in $M$ if and only if $g(h(D^\perp,D),\phi D^\perp)=0.$
\\\
\\\
\textbf{Theorem 2.2.6.} \textit{Let $M$ be a $\xi$-horizontal CR-submanifold of an LP-cosymplectic manifold $\tilde{M}$. Then $M$ is a contact CR-product if and only if $A_{\phi W}X=0~~\forall X\in D,W\in D^\perp$.}\\\
\\\
\textit{Proof.} Let $M$ be a contact CR-product. From Theorem 2.2.5 we have $\forall W,Z\in D^\perp,Y\in D,$\\\
\centerline{$0=g(h(W,Y),\phi Z)=g(A_{\phi Z}Y,W)$.}\\\
\par{Interchanging $Z,W$ and replacing $Y$ by $X$ we get,}
\begin{equation}
0=g(A_{\phi W}X,Z)\Rightarrow A_{\phi W}X\in D~~~\forall X\in D,W\in D^\perp. \tag{2.29}
\end{equation}
\par{As the leaf $M^T$ of $D$ is totally geodesic in $M$, for $Y\in D$,}\\\
$g(A_{\phi W}X,Y)=g(h(X,Y),\phi W)=g(\phi h(X,Y),W)=g(\phi\bar{\nabla}_XY,W)\\=g(\bar{\nabla}_X\phi Y,W)=g(\nabla_X\phi Y,W)=0$\\\
$\Rightarrow A_{\phi W}X\in D^\perp$. $\qquad \qquad \qquad \qquad \qquad \qquad \qquad \qquad \qquad \qquad \qquad $ (2.30)\\\
\par{From (2.29) and (2.30) we have, $A_{\phi W}X=0$.}\\\
\par{Conversely, $\forall X,Y\in D$ and $\forall Z,W\in D^\perp$,\\\
$0=g(A_{\phi W}X,Z)=g(h(X,Z),\phi W)$}\\\
$\Rightarrow$leaf $M^\perp$ of $D^\perp$ is totally geodesic in $M$ by Theorem 2.2.5.\\\
$g(\nabla_XY,Z)=g(\bar{\nabla}_XY,Z)=g(\phi\bar{\nabla}_XY,\phi Z)=g(\bar{\nabla}_X\phi Y,\phi Z)=g(\nabla_X\phi Y\\+h(X,\phi Y),\phi Z)=g(h(X,\phi Y),\phi Z)=g(A_{\phi Z}X,\phi Y)=0$\\\
$\Rightarrow$leaf $M^T$ of $D$ is totally geodesic in $M$. Thus the submanifold $M$ is a contact CR-product.\\\
\
\subsection{\large CR-submanifolds of an S-manifold}
In this section I have stated and proved some results related to the CR-submanifolds of an S-manifold.\\\
\\\
\\\
\textbf{Theorem 2.3.1.} \textit{Let $M$ be a totally umbilical contact CR-submanifold of an S-manifold $\tilde{M}$. Then the anti-invariant distribution $D^\perp$ is 1-dimensional, i.e., dim$D^\perp$=1.}\\\
\\\
\textit{Proof.} From (1.24) we have $\forall Z,W\in D^\perp$,\\\
$(\bar{\nabla}_Z\phi)W=\sum_{\alpha}\{g(\phi Z,\phi W)\xi_\alpha+\eta_\alpha(W)\phi^2Z\}$\\\
$\Rightarrow\bar{\nabla}_Z\phi W=\phi\bar{\nabla}_ZW+\sum_{\alpha}\{g(\phi Z,\phi W)\xi_\alpha+\eta_\alpha(W)\phi^2Z\}$.\\\
\par{Since $M$ is totally umbilical, by Weingarten and Gauss formula we have,}\\\
$\nabla^\perp_Z\phi W-g(H,\phi W)Z=\phi[\nabla_ZW+g(Z,W)H]+\sum_{\alpha}\{g(\phi Z,\phi W)\xi_\alpha\\+\eta_\alpha(W)\phi^2Z\}$.\\\
\par{Taking inner product with $Z$ we get,}\\\
$-g(H,\phi W)\|Z\|^2=g(Z,W)g(\phi H,Z)+\sum_{\alpha}\{g(\phi Z,\phi W)\eta_\alpha(Z)\\+\eta_\alpha(W)g(\phi^2Z,Z)\}$.$\qquad \qquad \qquad \qquad \qquad \qquad \qquad \qquad \qquad \qquad \qquad$(2.31)\\\
\par{Interchanging $Z,W$ in (2.31) and using (2.31) on the obtained equation we get after some calculations,}\\\
$g(H,\phi Z)\Big[\|W\|^2-\frac{g(Z,W)^2}{\|Z\|^2}\Big]+\frac{g(Z,W)}{\|Z\|^2}[-\sum_{\alpha}\eta_\alpha(W)\|\phi Z\|^2\\+\sum_{\alpha}g(\phi Z,\phi W)\eta_{\alpha}(Z)]+\sum_{\alpha}g(\phi W,\phi Z)\eta_{\alpha}(W)-\sum_{\alpha}\eta_\alpha(Z)\|\phi W\|^2=0$\\\
...it has a solution if $Z||W$, then dim$D^\perp$=1.\\\
\\\
\\\
\textbf{Corollary 2.3.1.} \textit{A contact CR-submanifold of an S-manifold reduces to a cosymplectic manifold provided the vector field $Z$ becomes a structure vector field $\xi_\beta$.}\\\
\\\
\textit{Proof.} From (1.24) we have,\\\
$(\bar{\nabla}_Z\phi)W=\sum_{\alpha}\{g(\phi Z,\phi W)\xi_{\alpha}+\eta_{\alpha}(W)\phi^2Z\}$.\\\
\par{Taking inner product with $Z$ we get,}\\\
$g((\bar{\nabla}_Z\phi)W,Z)=g(\phi Z,\phi W)\sum_{\alpha}\eta_{\alpha}(Z)+g(\phi^2Z,Z)\sum_{\alpha}\eta_{\alpha}(W)$\\\
$\Rightarrow g((\bar{\nabla}_\xi\phi)W,\xi_\beta)=g(\phi\xi_\beta,\phi W)\sum_{\alpha}\eta_\alpha(\xi_\beta)+g(\phi^2\xi_\beta,\xi_\beta)\sum_{\alpha}\eta_\alpha(W)$\\\
$~~~~~~~~~~~~~~~~~~~~~~~~~$ $=0$ (since $\phi\xi_\beta=0$).\\\
\par{Hence the manifold becomes cosymplectic.}\\\
\\\
\\\
\textbf{Theorem 2.3.2.} \textit{Let $M$ be a CR-submanifold of an S-manifold $\tilde{M}$. Then the leaf $M^\perp$ of $D^\perp$ is totally geodesic in $M$ if and only if $g(h(W,Y),QZ)\\+g(\phi W,\phi Z)\sum_{\alpha}\eta_\alpha(Y)=0~~\forall Y\in D$ and $\forall Z,W\in D^\perp$.}\\\
\\\
\textit{Proof.} Using (2.14) and (1.24) we get,\\\
$(\bar{\nabla}_XP)Y+(\bar{\nabla}_XQ)Y+h(X,PY)-A_{QY}X-Bh(X,Y)-Ch(X,Y)$\\\
$=g(\phi X,\phi Y)\sum_{\alpha}\xi_{\alpha}+\phi^2X\sum_{\alpha}\eta_{\alpha}(Y)$\\\
$\Rightarrow(\bar{\nabla}_XP)Y=A_{QY}X+Bh(X,Y)+g(\phi X,\phi Y)\sum_{\alpha}\xi_{\alpha}+\phi PX\sum_{\alpha}\eta_{\alpha}(Y)$,\\\
and $(\bar{\nabla}_XQ)Y=Ch(X,Y)-h(X,PY)+\phi QX\sum_{\alpha}\eta_{\alpha}(Y)$.\\\
\par{Thus, $P\nabla_XY=\nabla_XPY-A_{QY}X-Bh(X,Y)-g(\phi X,\phi Y)\sum_{\alpha}\xi_{\alpha}\\-\phi PX\sum_{\alpha}\eta_{\alpha}(Y)$,$\qquad \qquad \qquad \qquad \qquad \qquad \qquad \qquad \qquad \qquad \qquad$(2.32)\\\
and $Q\nabla_XY=\nabla^\perp_XQY+h(X,PY)-\phi QX\sum_{\alpha}\eta_{\alpha}(Y)-Ch(X,Y)$}.\\\
\par{Putting $X=W,Y=Z\in D^\perp$ in (2.32) we get,}\\\
$P\nabla_WZ=-A_{QZ}W-Bh(W,Z)-g(\phi W,\phi Z)\sum_{\alpha}\xi_{\alpha}$.\\\
\par{Taking inner product with $Y\in D$ and using (1.34) we have,}\\\
$g(\nabla_WZ,PY)=-g(A_{QZ}W,Y)-g(\phi W,\phi Z)\sum_{\alpha}\eta_{\alpha}(Y)$\\\
$~~~~~~~~~~~~~~~~~~$ $=-g(h(W,Y),QZ)-g(\phi W,\phi Z)\sum_{\alpha}\eta_{\alpha}(Y)$\\\
$\Rightarrow$the leaf $M^\perp$ of $D^\perp$ is totally geodesic in $M$ if and only if $g(h(W,Y),QZ)\\+g(\phi W,\phi Z)\sum_{\alpha}\eta_\alpha(Y)=0~~\forall Y\in D$ and $\forall Z,W\in D^\perp$.\\\
\
\subsection{\large CR-submanifolds of a GKM-manifold}
In this section I have stated and proved two theorems related to the structure vector field and leaf respectively regarding a CR-submanifold of a GKM-manifold. At last, I have given an example of a GKM manifold.\\\
\\\
\\\
\textbf{Theorem 2.4.1.} \textit{A contact CR-submanifold of a GKM-manifold reduces to a cosymplectic manifold provided the vector field $Z$ becomes the structure vector field $\xi_j$.}\\\
\\\
\textit{Proof.} From (1.26) we have $\forall Z,W\in D^\perp$,\\\
\centerline{$(\bar{\nabla}_Z\phi)W=\sum_{i=1}^{s}\{g(\phi Z,W)\xi_i-\eta_i(W)\phi^2Z\}$.}\\\
\par{Taking inner product with $Z$ we get,}\\\
\centerline{$g((\bar{\nabla}_Z\phi)W,Z)=g(\phi Z,W)\eta(Z)-g(\phi^2Z,Z)\sum_{i=1}^{s}\eta_i(W)$.}\\\
\par{Putting $Z=\xi_j$ we have,}\\\
\centerline{$g((\bar{\nabla}_{\xi_j}\phi)W,\xi_j)=g(\phi\xi_j,W)\eta(\xi_j)-g(\phi^2\xi_j,\xi_j)\sum_{i=1}^{s}\eta_i(W)=0$}\\\
(since $\phi\xi_j=0$).\\\
\par{Hence the proof.}\\\
\\\
\\\
\textbf{Theorem 2.4.2.} \textit{Let $M$ be a CR-submanifold of a GKM-manifold $\tilde{M}$. Then the leaf $M^\perp$ of $D^\perp$ is totally geodesic in $M$ if and only if}\\\
$g(h(W,Y),QZ)+g(\phi W,Z)\sum_{i=1}^{s}\eta_i(Y)=0~~~\forall Y\in D$ and $\forall Z,W\in D^\perp$.\\\
\\\
\textit{Proof.} Using (2.14) and (1.26) we get,\\\
$(\bar{\nabla}_XP)Y+(\bar{\nabla}_XQ)Y+h(X,PY)-A_{QY}X-Bh(X,Y)-Ch(X,Y)$\\\
$=\sum_{i=1}^{s}\{g(\phi X,Y)\xi_i-\eta_i(Y)\phi^2X\}$\\\
$\Rightarrow(\bar{\nabla}_XP)Y=A_{QY}X+Bh(X,Y)+g(\phi X,Y)\sum_{i=1}^{s}\xi_i-\phi PX\sum_{i=1}^{s}\eta_i(Y)$,\\\
and $(\bar{\nabla}_XQ)Y=Ch(X,Y)-h(X,PY)-\phi QX\sum_{i=1}^{s}\eta_i(Y)$.\\\
\par{Thus $P\nabla_XY=\nabla_XPY-A_{QY}X-Bh(X,Y)-g(\phi X,Y)\sum_{i=1}^{s}\xi_i\\+\phi PX\sum_{i=1}^{s}\eta_i(Y)$.}\\\
\par{Hence $P\nabla_WZ=-A_{QZ}W-Bh(W,Z)-g(\phi W,Z)\sum_{i=1}^{s}\xi_i~~\forall W,Z\in D^\perp$.}\\\
\par{Taking inner product with $Y\in D$ and using (1.34) we get,}\\\
$g(\nabla_WZ,PY)=-g(A_{QZ}W,Y)-g(\phi W,Z)\sum_{i=1}^{s}\eta_i(Y)$\\\
$~~~~~~~~~~~~~~~~$ $=-g(h(W,Y),QZ)-g(\phi W,Z)\sum_{i=1}^{s}\eta_i(Y)$\\\
$\Rightarrow$the leaf $M^\perp$ of $D^\perp$ is totally geodesic in $M$ if and only if
$g(h(W,Y),QZ)\\+g(\phi W,Z)\sum_{i=1}^{s}\eta_i(Y)=0~~\forall Y\in D$ and $\forall Z,W\in D^\perp$.
\\\
\\\
\\\
$\bullet$ \underline{\textbf{Example of a GKM-manifold}} \textbf{:}\\\
\
\par{We take $n=2,s=3$ and the followings$-$ }\\\
$e_1=f_1(z_1,z_2,z_3)\frac{\partial}{\partial x_1}+f_2(z_1,z_2,z_3)\frac{\partial}{\partial y_1}$,
$e_2=-f_2(z_1,z_2,z_3)\frac{\partial}{\partial x_2}+f_1(z_1,z_2,z_3)\frac{\partial}{\partial y_1}$,\\\
$e_3=f_1(z_1,z_2,z_3)\frac{\partial}{\partial x_2}+f_2(z_1,z_2,z_3)\frac{\partial}{\partial y_2}$,
$e_4=-f_2(z_1,z_2,z_3)\frac{\partial}{\partial x_2}+f_1(z_1,z_2,z_3)\frac{\partial}{\partial y_2}$,\\\
$e_5=\frac{\partial}{\partial z_1}$,
$e_6=\frac{\partial}{\partial z_2}$,
$e_7=\frac{\partial}{\partial z_3}$,\\\
$f_1=c_2e^{-(z_1+z_2+z_3)}cos(z_1+z_2+z_3)-c_1e^{-(z_1+z_2+z_3)}sin(z_1+z_2+z_3)$,\\\
$f_2=c_1e^{-(z_1+z_2+z_3)}cos(z_1+z_2+z_3)+c_2e^{-(z_1+z_2+z_3)}sin(z_1+z_2+z_3)$,\\\
$g=\frac{1}{f_1^2+f_2^2}\sum_{i=1}^{s}(dx_i\otimes dx_i+dy_i\otimes dy_i+dz_1\otimes dz_1)+dz_2\otimes dz_2+dz_3\otimes dz_3$,\\\
$\eta_1(X)=g(X,e_5),\eta_2(X)=g(X,e_6),\eta_3(X)=g(X,e_7)$ for any vector field $X$ on $M$,\\\
$\phi e_1=e_2,\phi e_2=-e_1,\phi e_3=e_4,\phi e_4=-e_3,\phi(e_5=\xi_1)=0,\phi(e_6=\xi_2)=0,\phi(e_7=\xi_3)=0$.\\\
\par{Here $M$ becomes a GKM-manifold with structure tensors $(\phi ,\xi_i ,\eta_i ,g)$, $i=1,2,3.$}\\\
\section{\large Acknowledgement}
The author has been sponsored by University Grants Commission(UGC) Junior Research Fellowship, India. \textit{UGC-Ref. No.}: 1139/(CSIR-UGC NET JUNE 2018).\\\

\end{document}